\documentclass[12pt]{article}
\usepackage{amsmath,amsfonts,amssymb,amsthm}
\usepackage{fullpage}

\title{The local index formula for $SU_q(2)$}
\author{Walter van Suijlekom,$^1$ \
Ludwik D\c{a}browski,$^1$ \
Giovanni Landi,$^2$ \\
Andrzej Sitarz,$^3$\footnotemark\ \
Joseph C. V\'arilly$^4$\footnotemark \\[12pt]
$^1\,$Scuola Internazionale Superiore di Studi Avanzati,\\
      Via Beirut 2--4, 34014 Trieste, Italy\\[3pt]
$^2\,$Dipartimento di Matematica e Informatica,
      Universit\`a di Trieste,\\
      Via Valerio 12/1, 34127 Trieste\\
and INFN, Sezione di Napoli, Napoli, Italy\\[3pt]
$^3\,$Institute of Physics, Jagiellonian University,\\
      Reymonta 4, 30--059 Krak\'ow, Poland\\[3pt]
$^4\,$Departamento de Matem\'atica,
      Universidad de Costa Rica,\\
      2060 San Jos\'e, Costa Rica}

\date{}

\makeatletter
\def\section{\@startsection{section}{1}{\z@}{-3.5ex plus -1ex minus
  -.2ex}{2.3ex plus .2ex}{\large\bf}}
\def\subsection{\@startsection{subsection}{2}{\z@}{-3.25ex plus -1ex
  minus -.2ex}{1.5ex plus .2ex}{\normalsize\bf}}
\makeatother

\newtheorem{thm}{Theorem}[section]
\newtheorem{prop}[thm]{Proposition}
\newtheorem{lem}[thm]{Lemma}

\theoremstyle{definition}
\newtheorem{defn}[thm]{Definition}

\newtheorem{rem}{Remark}
\numberwithin{equation}{section}

\newcommand{\A}{\mathcal{A}}        
\renewcommand{\a}{\alpha}           
\newcommand{\B}{\mathcal{B}}        
\newcommand{\braket}[2]{\langle#1\mathbin|#2\rangle} 
\newcommand{\C}{\mathbb{C}}         
\newcommand{\cha}[1]{C^{#1}(\mathcal{A})}
\newcommand{\chala}[1]{C^{#1}_\lambda(\mathcal{A})}

\newcommand{\dn}{{\mathord{\downarrow}}} 
\DeclareMathOperator{\Dom}{Dom}     
\newcommand{\eps}{\varepsilon}      
\newcommand{\ev}{{\mathrm{ev}}}     
\newcommand{\ga}{\gamma}            
\renewcommand{\H}{\mathcal{H}}      
\newcommand{\half}{{\mathchoice{\oh}{\oh}{\shalf}{\shalf}}} 
\DeclareMathOperator{\ind}{Index}   
\newcommand{\K}{\mathcal{K}}        
\newcommand{\ket}[1]{|#1\rangle}    
\newcommand{\kett}[1]{|#1\rangle\!\rangle} 
\renewcommand{\L}{\mathcal{L}}      
\newcommand{\la}{\lambda}           
\DeclareMathOperator{\Mat}{Mat}     
\newcommand{\N}{\mathbb{N}}         
\newcommand{\nn}{\nonumber}         
\newcommand{\odd}{{\mathrm{odd}}}   
\newcommand{\oh}{{\tfrac{1}{2}}}    
\newcommand{\ooh}{{\tfrac{3}{2}}}   
\newcommand{\OP}{\mathrm{OP}}       
\newcommand{\otto}{\leftrightarrow} 
\newcommand{\ox}{\otimes}           
\newcommand{\piappr}{\underline{\pi}} 
\newcommand{\piso}[1]{\lfloor#1\rfloor} 
\newcommand{\R}{\mathbb{R}}         
\DeclareMathOperator{\Res}{Res}     
\newcommand{\sepword}[1]{\quad\mbox{#1}\quad} 
\newcommand{\set}[1]{\{\,#1\,\}}    
\newcommand{\Sf}{\mathbb{S}}        
\newcommand{\sg}{\sigma}            
\newcommand{\shalf}{{\scriptstyle\frac{1}{2}}} 
\DeclareMathOperator{\Sign}{Sign}   
\newcommand{\ssesq}{{\scriptstyle\frac{3}{2}}} 
\newcommand{\sesq}{{\mathchoice{\ooh}{\ooh}{\ssesq}{\ssesq}}} 
\DeclareMathOperator{\Tr}{Tr}       
\newcommand{\ul}[1]{\underline{#1}} 
\newcommand{\up}{{\mathord{\uparrow}}} 
\newcommand{\x}{\times}             
\newcommand{\Z}{\mathbb{Z}}         
\newcommand{\8}{\bullet}            
\renewcommand{\:}{\colon}           
\def\<#1,#2>{\langle#1,#2\rangle}   

\newbox\ncintdbox \newbox\ncinttbox
\setbox0=\hbox{$-$}
\setbox2=\hbox{$\displaystyle\int$}
\setbox\ncintdbox=\hbox{\rlap{\hbox
    to \wd2{\kern-.1em\box2\relax\hfil}}\box0\kern.1em}
\setbox0=\hbox{$\vcenter{\hrule width 4pt}$}
\setbox2=\hbox{$\textstyle\int$}
\setbox\ncinttbox=\hbox{\rlap{\hbox
    to \wd2{\kern-.05em\box2\relax\hfil}}\box0\kern.1em}
\newcommand{\ncint}{\mathop{\mathchoice{\copy\ncintdbox}%
    {\copy\ncinttbox}{\copy\ncinttbox}%
    {\copy\ncinttbox}}\nolimits}

\begin{document}

\maketitle

\thispagestyle{empty}

\begin{abstract}
We discuss the local index formula of Connes--Moscovici for the
isospectral noncommutative geometry that we have recently constructed
on quantum $SU(2)$. We work out the cosphere bundle and the dimension
spectrum as well as the local cyclic cocycles yielding the index
formula.
\end{abstract}

\vfill

\textit{Key words and phrases}:
Noncommutative geometry, spectral triple, quantum $SU(2)$.

\textit{Mathematics Subject Classification:}
Primary 58B34; Secondary 17B37.

\renewcommand{\thefootnote}{\fnsymbol{footnote}}
\addtocounter{footnote}{1}
\footnotetext{Partially supported by Polish State
Committee for Scientific Research (KBN) under grant 2\,P03B\,022\,25 and by PBZ-MIN-008/P03/2003 grant.}
\addtocounter{footnote}{1}
\footnotetext{Regular Associate of the Abdus Salam ICTP, Trieste.}

\newpage

\section{Introduction}

Recent investigations show that the ``quantum space'' underlying the
quantum group $SU_q(2)$ is an important arena for testing and
implementing ideas coming from noncommutative differential geometry.
In \cite{Naiad} it has been endowed with an isospectral tridimensional
geometry via a bi-equivariant $3^+$-summable spectral triple $(\A(SU_q(2)),\H,D)$.
Earlier, a ``singular'' (in the sense of not admitting a commutative
limit) spectral triple was constructed in~\cite{ChakrabortyPEqvt}. The
latter geometry was put in the general theory of Connes--Moscovici
\cite{ConnesMIndex} by a systematic discussion of the local index
formula \cite{ConnesSUq}. In this paper, we present a similar analysis
for the former geometry. It turns out that most of the results
coincide with those of~\cite{ConnesSUq}.

The main idea of that paper is to construct a (quantum) cosphere
bundle $\Sf_q^*$ on $SU_q(2)$, that considerably simplifies the
computations concerning the local index formula. Essentially, with the
operator derivation $\delta$ defined by
$\delta(T) := |D|T - T|D|$, one considers an operator $x$ in the
algebra $\B = \bigcup_{n=0}^\infty \delta^n(\A)$ up to smoothing
operators; these give no contribution to the residues appearing in the
local cyclic cocycle giving the local index formula. The removal of
the irrelevant smoothing operators is accomplished by introducing a
symbol map from $SU_q(2)$ to the cosphere bundle $\Sf_q^*$. The latter
is defined by its algebra $C^\infty(\Sf_q^*)$ of ``smooth functions''
which is, by definition, the image of a map
$$
\rho : \B \to C^\infty(D^2_{q+} \x D^2_{q-} \x \Sf^1)
$$
where $D^2_{q\pm}$ are two quantum disks. One finds that an element
$x$ in the algebra $\B$ can be determined up to smoothing operators
by~$\rho(x)$.

In our present case, the cosphere bundle coincides with the one
obtained in \cite{ConnesSUq}; the same being true for the dimension
spectrum. Indeed, using this much simpler form of operators up to
smoothing ones, it is not difficult to compute the dimension spectrum
and obtain simple expressions for the residues appearing in the local
index formula. We find that the dimension spectrum is simple and given
by the set $\{1,2,3\}$.

The cyclic cohomology of the algebra $\A(SU_q(2)$ has been computed
explicitly in \cite{MasudaNW} where it was found to be given in terms of a
single generator. We express this element in terms of a single local
cocycle similarly to the computations in \cite{ConnesSUq}. But
contrary to the latter, we get an extra term involving $P |D|^{-3}$
which drops in \cite{ConnesSUq}, being traceclass for the case
considered there. Here $P = \half(1 + F)$ with $F = \Sign D$, the sign
of the operator~$D$.

Finally as a simple example, we compute the Fredholm index of $D$
coupled with the unitary representative of the generator of
$K_1(\A(SU_q(2)))$.

\section{The isospectral geometry of \boldmath $SU_q(2)$}
\label{sec:iso-ge}

We recall the construction of the spectral triple $(\A(SU_q(2)),\H,D)$
of \cite{Naiad}. Let $\A = \A(SU_q(2))$ be the $*$-algebra generated
by $a$ and~$b$, subject to the following commutation rules:
\begin{gather}
ba = q ab,  \qquad  b^*a = qab^*, \qquad bb^* = b^*b,
\nn \\
a^*a + q^2 b^*b = 1,  \qquad  aa^* + bb^* = 1.
\label{eq:suq2-relns}
\end{gather}
In the following we shall take $0 < q < 1$. Note that we have
exchanged $a \otto a^*$, $b \otto -b$ with respect to the notation of
\cite{ChakrabortyPEqvt} and~\cite{ConnesSUq}.

The Hilbert space of spinors $\H$ has an orthonormal basis labelled as
follows. For each $j = 0,\half,1,\dots$, we abbreviate
$j^+ = j + \half$ and $j^- = j - \half$. The orthonormal basis
consists of vectors $\ket{j\mu n\up}$ for $j = 0,\half,1,\dots$,
$\mu = -j,\dots,j$ and $n = -j^+,\dots,j^+$; together with
$\ket{j\mu n\dn}$ for $j=\half,1,\dots$, $\mu = -j,\ldots,j$ and
$n = -j^-,\dots,j^-$. We adopt a vector notation by juxtaposing the
pair of spinors
\begin{equation}
\kett{j\mu n} := \begin{pmatrix} \ket{j\mu n\up} \\[2\jot]
\ket{j\mu n\dn} \end{pmatrix},
\label{eq:kett-defn}
\end{equation}
and with the convention that the lower component is zero when
$n = \pm(j + \half)$ or $j = 0$. In this way, we get a decomposition
$\H = \H^\up \oplus \H^\dn$ into subspaces spanned by the ``up'' and
``down'' kets, respectively.

The spinor representation is the $*$-representation $\pi$ of $\A$ on
$\H$ --denoted by $\pi'$ in~\cite{Naiad}-- defined as follows. We set
$\pi(a) := a_+ + a_-$ and $\pi(b) := b_+ + b_-$, where $a_\pm$ and
$b_\pm$ are the following operators in $\H$:
\begin{align}
a_+ \,\kett{j\mu n}
&:= q^{(\mu+n-\half)/2} [j + \mu + 1]^\half
    \begin{pmatrix}
    q^{-j-\half} \, \frac{[j+n+\sesq]^{1/2}}{[2j+2]} & 0 \\[2\jot]
    q^\half \,\frac{[j-n+\half]^{1/2}}{[2j+1]\,[2j+2]} &
    q^{-j} \, \frac{[j+n+\half]^{1/2}}{[2j+1]}
    \end{pmatrix} \kett{j^+ \mu^+ n^+},
\nn \\
a_- \,\kett{j\mu n}
&:= q^{(\mu+n-\half)/2} [j - \mu]^\half
    \begin{pmatrix}
    q^{j+1} \, \frac{[j-n+\half]^{1/2}}{[2j+1]} &
    - q^\half \,\frac{[j+n+\half]^{1/2}}{[2j]\,[2j+1]} \\[2\jot]
    0 & q^{j+\half} \, \frac{[j-n-\half]^{1/2}}{[2j]}
    \end{pmatrix} \kett{j^- \mu^+ n^+},
\nn \\
b_+ \,\kett{j\mu n}
&:= q^{(\mu+n-\half)/2} [j + \mu + 1]^\half
    \begin{pmatrix}
    \frac{[j-n+\sesq]^{1/2}}{[2j+2]} & 0 \\[2\jot]
    - q^{-j-1} \,\frac{[j+n+\half]^{1/2}}{[2j+1]\,[2j+2]} &
    q^{-\half} \, \frac{[j-n+\half]^{1/2}}{[2j+1]}
    \end{pmatrix} \kett{j^+ \mu^+ n^-},
\nn \\
b_- \,\kett{j\mu n}
&:= q^{(\mu+n-\half)/2} [j - \mu]^\half
    \begin{pmatrix}
    - q^{-\half} \, \frac{[j+n+\half]^{1/2}}{[2j+1]} &
    - q^j \,\frac{[j-n+\half]^{1/2}}{[2j]\,[2j+1]} \\[2\jot]
    0 & - \frac{[j+n-\half]^{1/2}}{[2j]}
    \end{pmatrix} \kett{j^- \mu^+ n^-}.
\label{eq:spin-rep}
\end{align}
Here $[N] := (q^{-N} - q^N)/(q^{-1} - q)$ is a ``$q$-integer''.

The Dirac operator $D$ that was exhibited in \cite{Naiad} is diagonal
in the given orthonormal basis of~$\H$, and is one of a family of
selfadjoint operators of the form
\begin{equation}
D \kett{j\mu n} = \begin{pmatrix} d^\up j + c^\up & 0 \\
0 & d^\dn j + c^\dn \end{pmatrix} \kett{j\mu n},
\label{eq:linear-evs}
\end{equation}
where $d^\up, d^\dn, c^\up, c^\dn$ are real numbers not depending 
on~$j,\mu,n$. In order that the sign of~$D$ be nontrivial we need to
assume $d^\dn d^\up < 0$, so we may as well take $d^\up > 0$ and
$d^\dn < 0$.

Apart from the issue of their signs, the particular constants that
appear in \eqref{eq:linear-evs} are fairly immaterial: $c^\up$ and
$c^\dn$ do not affect the index calculations later on while $d^\up$ and
$|d^\dn|$ yield scaling factors on some noncommutative integrals. Thus
little generality is lost by making the following choice,
\begin{equation}
D \kett{j\mu n} = \begin{pmatrix} 2j + \sesq & 0 \\
0 & -2j - \half \end{pmatrix} \kett{j\mu n}.
\label{eq:classical-evs}
\end{equation}
whose spectrum (with multiplicity!) coincides with that of the
classical Dirac operator of the sphere $\Sf^3$ equipped with the
round metric (indeed, the spin geometry of the 3-sphere can now be
recovered by taking $q = 1$).

We let $D = F\,|D|$ be the polar decomposition of $D$ where
$|D| := (D^2)^\half$ and $F = \Sign D$. Explicitly, we see that
$$
F\kett{j\mu n}
= \begin{pmatrix} 1 & 0 \\ 0 & -1 \end{pmatrix} \kett{j\mu n},
\qquad
|D|\,\kett{j\mu n} = \begin{pmatrix} 2j + \sesq & 0 \\
0 & 2j + \half \end{pmatrix} \kett{j\mu n}.
$$
Clearly, $P^\up := \half(1 + F)$ and
$P^\dn := \half(1 - F) = 1 - P^\up$ are the orthogonal projectors
whose range spaces are $\H^\up$ and $\H^\dn$, respectively.

\begin{prop}
\label{pr:regular}
The triple $(\A(SU_q(2)),\H,D)$ is a regular $3^+$-summable
spectral triple.
\end{prop}

\begin{proof}
It was already shown in \cite{Naiad} that this spectral triple is
$3^+$-summable: indeed, this follows easily from the growth of the
eigenvalues in~\eqref{eq:classical-evs}. The remaining issue is its
regularity. Recall \cite{CareyPRS,ConnesMIndex,Polaris} that this
means that the algebra generated by $\A$ and $[D,\A]$ should lie
within the smooth domain $\bigcap_{n=0}^\infty \Dom \delta^n$ of the
operator derivation $\delta(T) := |D|T - T|D|$.

Since $2j + \sesq = 2j^+ + \half$ and $2j + \half = 2j^- + \sesq$ and
due to the triangular forms of the matrices in \eqref{eq:spin-rep},
the off-diagonal terms vanish in the $2 \x 2$-matrix expressions for
$\delta(a_+)$ and $\delta(a_-)$. Indeed one finds,
\begin{align*}
\delta(a_+)\kett{j\mu n}
&= \begin{pmatrix} 2j + \tfrac{5}{2} & 0 \\
0 & 2j + \sesq \end{pmatrix} a_+ \kett{j\mu n}
- a_+ \begin{pmatrix} 2j + \sesq & 0 \\
0 & 2j + \half \end{pmatrix} \kett{j\mu n},
\\
\delta(a_-)\kett{j\mu n}
&= \begin{pmatrix} 2j + \half & 0 \\
0 & 2j - \half \end{pmatrix} a_- \kett{j\mu n}
- a_- \begin{pmatrix} 2j + \sesq & 0 \\
0 & 2j + \half \end{pmatrix} \kett{j\mu n}.
\end{align*}
In both cases we obtain
\begin{equation}
\delta(a_+) = P^\up a_+ P^\up + P^\dn a_+ P^\dn, \qquad
\delta(a_-) = - P^\up a_- P^\up - P^\dn a_- P^\dn.
\label{eq:delta-apm}
\end{equation}
Replacing $a$ by $b$, the same triangular matrix structure leads to
\begin{equation}
\delta(b_+) = P^\up b_+ P^\up + P^\dn b_+ P^\dn, \qquad
\delta(b_-) = - P^\up b_- P^\up - P^\dn b_- P^\dn.
\label{eq:delta-bpm}
\end{equation}
Thus $\delta(\pi(a)) = \delta(a_+) + \delta(a_-)$ is bounded, with
$\|\delta(\pi(a))\| \leq \|\pi(a)\|$; and likewise for $\pi(b)$.
Next, $\delta([D,a_+]) = [D,\delta(a_+)]$, so that
$$
\delta([D,a_+]) \kett{j\mu n}
= \begin{pmatrix} 2j + \tfrac{5}{2} & 0 \\
0 & -2j - \sesq \end{pmatrix} \delta(a_+) \kett{j\mu n}
- \delta(a_+) \begin{pmatrix} 2j + \sesq & 0 \\
0 & -2j - \half \end{pmatrix} \kett{j\mu n},
$$
since all matrices appearing are diagonal. This, together with the
analogous calculation for $\delta([D,a_-])$, shows that
\begin{equation}
\delta([D,a_+]) = P^\up a_+ P^\up - P^\dn a_+ P^\dn,  \qquad
\delta([D,a_-]) = P^\up a_- P^\up - P^\dn a_- P^\dn.
\label{eq:delta-Dapm}
\end{equation}
A similar argument for $b$ gives
\begin{equation}
\delta([D,b_+]) = P^\up b_+ P^\up - P^\dn b_+ P^\dn,  \qquad
\delta([D,b_-]) = P^\up b_- P^\up - P^\dn b_- P^\dn.
\label{eq:delta-Dbpm}
\end{equation}
Combining \eqref{eq:delta-apm}, \eqref{eq:delta-Dapm}, and the
analogous relations with $a$ replaced by~$b$, we see that both
$\A$ and $[D,\A]$ lie within $\Dom \delta$. An easy induction
shows that they also lie within $\Dom \delta^k$ for $k = 2,3,\dots$.
\end{proof}

This proposition continues to hold if we replace $\A(SU_q(2))$ by a
suitably completed algebra, which is stable under the holomorphic
function calculus.

\vspace{6pt}

Let $\Psi^0(\A)$ be the algebra generated by $\delta^k(\A)$ and
$\delta^k([D,\A])$ for all $k \geq 0$ (the notation suggests that, in
the spirit of~\cite{ConnesMIndex} one thinks of it as an ``algebra of
pseudodifferential operators of order~$0$''). Since, for
instance,
\begin{align*}
P^\up \pi(a) P^\up
&= \half \delta^2(\pi(a)) + \half \delta([D,\pi(a)]),
\\
P^\up a_+ P^\up
&= \half P^\up \pi(a) P^\up + \half P^\up \delta(\pi(a)) P^\up,
\end{align*}
we see that $\Psi^0(\A)$ is in fact generated by the diagonal-corner
operators $P^\up a_\pm P^\up$, $P^\dn a_\pm P^\dn$,
$P^\up b_\pm P^\up$, $P^\dn b_\pm P^\dn$ together with the
other-corner operators $P^\dn a_+ P^\up$, $P^\up a_- P^\dn$,
$P^\dn b_+ P^\up$, and $P^\up b_- P^\dn$.
Following~\cite{ConnesSUq}, let $\B$ be the algebra generated by all
$\delta^n(\A)$ for $n \geq 0$. It is a subalgebra of $\Psi^0(\A)$ and it
is generated by the diagonal operators
\begin{equation}
\tilde a_\pm := \pm \delta(a_\pm)
= P^\up a_\pm P^\up + P^\dn a_\pm P^\dn,  \qquad
\tilde b_\pm := \pm \delta(b_\pm)
= P^\up b_\pm P^\up + P^\dn b_\pm P^\dn,
\label{eq:extended-pi}
\end{equation}
and by the off-diagonal operators $P^\dn a_+ P^\up + P^\up a_- P^\dn$
and $P^\dn b_+ P^\up + P^\up b_- P^\dn$. \\

For later convenience we shall introduce an approximate representation
$\piappr$ found in \cite{Naiad}, which coincides with $\pi$ up to
compact operators. Note first, that the off-diagonal coefficients in
\eqref{eq:spin-rep} give rise to smoothing operators in
$\OP^{-\infty}$ (see Appendix~\ref{sec:pdc}), due to the terms
appearing in their denominators; we can furthermore simplify the
diagonal terms.

We set $\piappr(a) := \ul{a}_+ + \ul{a}_-$ and
$\piappr(b) := \ul{b}_+ + \ul{b}_-$ with the following definitions:
\begin{align}
\ul{a}_+ \,\kett{j\mu n}
&:= \sqrt{1-q^{2j+2\mu+2}} \begin{pmatrix}
\sqrt{1- q^{2j+2n+3}} & 0 \\ 0 & \sqrt{1- q^{2j+2n+1}}
\end{pmatrix} \kett{j^+ \mu^+ n^+},
\nn \\
\ul{a}_- \,\kett{j\mu n}
&:= q^{2j+\mu+n+\half} \begin{pmatrix}
q & 0 \\ 0 & 1 \end{pmatrix} \kett{j^- \mu^+ n^+},
\nn \\
\ul{b}_+ \,\kett{j\mu n}
&:= q^{j+n-\half} \sqrt{1-q^{2j+2\mu+2}} \begin{pmatrix}
q & 0 \\ 0 & 1 \end{pmatrix} \kett{j^+ \mu^+ n^-},
\nn \\
\ul{b}_- \,\kett{j\mu n}
&:= -q^{j+\mu} \begin{pmatrix}
\sqrt{1-q^{2j+2n+1}} & 0 \\ 0 & \sqrt{1-q^{2j+2n-1}}
\end{pmatrix} \kett{j^- \mu^+ n^-}.
\label{eq:approx-repn}
\end{align}
These formulas can be obtained from \eqref{eq:spin-rep} by truncation,
using the pair of estimates
\begin{align*}
\bigl((q^{-1} - q)[n] \bigr)^{-1} - q^n &= q^{3n} + O(q^{5n}),
\\
1 - \sqrt{1 - q^\a} & \leq q^\a,  \qquad{\rm for ~any} \ \a \geq 0.
\end{align*}
The operators $\piappr(x) - \pi(x)$ are given by sequences of rapid
decay, and hence are elements in $\OP^{-\infty}$ (as defined in Appendix~\ref{sec:pdc}). Therefore, we can
replace $\pi$ by $\piappr$ when dealing with the local cocycle in the
local index theorem in the next section.

\begin{rem}
\label{rk:approx-repn}
These operators differ slightly from the approximate representation
given in \cite{Naiad}. Using the inequality
$1 - \sqrt{1-q^\alpha} \leq q^\alpha$, they can be seen to differ from
the operators therein by a compact operator in the principal ideal
$\K_q$ generated by the operator
$L_q \: \kett{j\mu n} \mapsto q^j \kett{j \mu n}$. Note that
$\K_q \subset \OP^{-\infty}$.
\end{rem}

Now, observe that
\begin{align}
[|D|, \piappr(a)] &= \ul{a}_+ - \ul{a}_-, &
[D, \piappr(a)] &= F(\ul{a}_+ - \ul{a}_-),
\nn \\
[|D|, \piappr(b)] &= \ul{b}_+ - \ul{b}_-, &
[D, \piappr(b)] &= F(\ul{b}_+ - \ul{b}_-),
\label{eq:comm-DA}
\end{align}
and also that $F$ commutes with $\ul{a}_\pm$ and $\ul{b}_\pm$.
The operators $\ul{a}_\pm$ and $\ul{b}_\pm$ have a simpler expression
if we use the following relabelling of the orthonormal basis of~$\H$,
\begin{align}
v_{xy\up}^j &:= \ket{j,x-j,y-j-\half,\up} \sepword{for}
x = 0, \dots, 2j; \ y = 0, \dots, 2j + 1,
\nn \\
v_{xy\dn}^j &:= \ket{j,x-j,y-j+\half,\dn} \sepword{for}
x = 0, \dots, 2j; \ y = 0, \dots, 2j - 1.
\label{eq:new-basis}
\end{align}
We again employ the pairs of vectors
$$
v_{xy}^j := \begin{pmatrix} v_{xy\up}^j \\ v_{xy\dn}^j \end{pmatrix},
$$
where the lower component is understood to be zero if $y = 2j$ or
$2j+1$, or if $j = 0$. The simplification is that on these vector
pairs, all the $2 \x 2$ matrices in \eqref{eq:approx-repn} become
scalar matrices,
\begin{align}
\ul{a}_+ v_{xy}^j
&= \sqrt{1 - q^{2x+2}} \sqrt{1 - q^{2y+2}} \, v_{x+1,y+1}^{j^+},
\nn \\
\ul{a}_- v_{xy}^j &= q^{x+y+1} \, v_{xy}^{j^-},
\nn \\
\ul{b}_+ v_{xy}^j &= q^y \sqrt{1 - q^{2x+2}} \, v_{x+1,y}^{j^+},
\nn \\
\ul{b}_- v_{xy}^j &= - q^x \sqrt{1 - q^{2y}} \, v_{x,y-1}^{j^-}.
\label{eq:approx-repn-bis}
\end{align}
These formulas coincide with those found in \cite[Sec.~6]{ConnesSUq}
up to a doubling of the Hilbert space and the change of conventions
$a \otto a^*$, $b \otto -b$. Indeed, since the spin representation is
isomorphic to a direct sum of two copies of the regular
representation, the formulas in \eqref{eq:approx-repn-bis} exhibit the
same phenomenon for the approximate representations.

\section{The cosphere bundle}
\label{sec:cosphere}

In \cite{ConnesSUq} Connes constructs a ``cosphere bundle'' using the
regular representation of $\A(SU_q(2))$. In view
of~\eqref{eq:approx-repn-bis}, the same cosphere bundle may be
obtained directly from the spin representation by adapting that
construction, as we now proceed to do. In what follows, we use the
algebra $\A = \A(SU_q(2))$, but we could as well replace it with its
completion $C^\infty(SU_q(2))$, which is closed under holomorphic
functional calculus (see Appendix~\ref{sec:pdc}).

We recall two well-known infinite dimensional representations
$\pi_\pm$ of $\A(SU_q(2))$ by bound\-ed operators on the Hilbert space
$\ell^2(\N)$. On the standard orthonormal basis
$\set{\eps_x : x \in \N}$, they are given by
\begin{equation}
\pi_\pm(a) \,\eps_x := \sqrt{1-q^{2x+2}} \,\eps_{x+1},  \qquad
\pi_\pm(b) \,\eps_x := \pm q^x \,\eps_x.
\label{eq:pi-pm}
\end{equation}
We may identify the Hilbert space $\H$ spanned by all $v_{xy\up}^j$
and $v_{xy\dn}^j$ with the subspace $\H'$ of
$\ell^2(\N)_x \ox \ell^2(\N)_y \ox \ell^2(\Z)_{2j} \ox \C^2$
determined by the parameter restrictions in~\eqref{eq:new-basis}.
Thereby, we get the correspondence
\begin{align}
\ul{a}_+ &\otto \pi_+(a) \ox \pi_-(a) \ox V \ox 1_2,
\nn \\
\ul{a}_- &\otto - q\,\pi_+(b) \ox \pi_-(b^*) \ox V^* \ox 1_2,
\nn \\
\ul{b}_+ &\otto - \pi_+(a) \ox \pi_-(b) \ox V \ox 1_2,
\nn \\
\ul{b}_- &\otto - \pi_+(b) \ox \pi_-(a^*) \ox V^* \ox 1_2,
\label{eq:corr-smoothing}
\end{align}
where $V$ is the unilateral shift operator
$\eps_{2j} \mapsto \eps_{2j+1}$ in $\ell^2(\Z)$. This again, apart from
the $2 \x 2$ identity matrix $1_2$, coincides with the formula (204)
in~\cite{ConnesSUq}, up to the aforementioned exchange of the
generators.

The shift $V$ in the action of the operators $\ul{a}_\pm$ and
$\ul{b}_\pm$ on $\H$ can be encoded using the $\Z$-grading coming
from the one-parameter group of automorphisms $\ga(t)$ generated
by~$|D|$,
\begin{equation}
\ga(t) = \begin{pmatrix} \ga_{\up\up}(t) & \ga_{\up\dn}(t) \\
\ga_{\dn\up}(t) & \ga_{\dn\dn}(t) \end{pmatrix},  \sepword{where}
\left\{ 
\begin{aligned}
\ga_{\up\up}(t)
&: P^\up T P^\up \mapsto P^\up e^{it|D|} T e^{-it|D|} P^\up,
\\
\ga_{\up\dn}(t)
&: P^\up T P^\dn \mapsto P^\up e^{it|D|} T e^{-it|D|} P^\dn,
\\
\ga_{\dn\up}(t)
&: P^\dn T P^\up \mapsto P^\dn e^{it|D|} T e^{-it|D|} P^\up,
\\
\ga_{\dn\dn}(t)
&: P^\dn T P^\dn \mapsto P^\dn e^{it|D|} T e^{-it|D|} P^\dn,
\end{aligned} \right.
\label{eq:spectral-flow}
\end{equation}
for any operator $T$ on $\H$. On the subalgebra of ``diagonal''
operators $T = P^\up T P^\up + P^\dn T P^\dn$, the compression
$\ga_{\up\up} \oplus \ga_{\dn\dn}$ detects the
shift of $j$ of the restrictions of~$T$ to $\H^\up$ and $\H^\dn$
respectively. For example, $\ga_{\up\up}(t) \oplus \ga_{\dn\dn}(t)
: a_\pm \mapsto e^{\pm it} a_\pm$, so that the $\Z$-grading encodes
the correct shifts $j \to j\pm \half$ in the formulas for $a_\pm$; and
likewise for $b_\pm$.

{}From equation \eqref{eq:pi-pm} it follows that
$b - b^* \in \ker \pi_\pm$, and so the representations $\pi_\pm$ are
not faithful on $\A(SU_q(2))$. We define two algebras
$\A(D^2_{q\pm})$ to be the corresponding quotients,
\begin{equation}
0 \to \ker \pi_\pm \to \A(SU_q(2))
\xrightarrow{r_\pm} \A(D^2_{q\pm}) \to 0.
\label{eq:q-disks}
\end{equation}
We elaborate a little on the structure of the algebras
$\A(D^2_{q\pm})$. For convenience, we shall omit the quotient maps
$r_\pm$ in this discussion. Then $b = b^*$ in $\A(D_{q\pm}^2)$, and
from the defining relations \eqref{eq:suq2-relns} of $\A(SU_q(2))$, we
obtain
\begin{gather}
b a = q\, a b,  \qquad    a^*b = q\, b a^*,
\nn \\
a^*a + q^2 b^2 = 1,  \qquad  aa^* + b^2 = 1.
\label{eq:Sq2-relns}
\end{gather}
These algebraic relations define two isomorphic quantum $2$-spheres
$\Sf^2_{q+} \simeq \Sf^2_{q-} =: \Sf^2_q$ which have a classical
subspace $\Sf^1$ given by the characters $b \mapsto 0$,
$a \mapsto \lambda$ with $|\lambda| = 1$. A substitution
$q \mapsto q^2$, followed by $b \mapsto q^{-2} b$ shows that $\Sf^2_q$
is none other than the equatorial Podle\'s sphere~\cite{Podles}. Thus,
the above quotients of $\A(SU_q(2))$ with respect to $\ker \pi_\pm$
either coincide with $\A(\Sf^2_q)$ or are quotients of it. Now, from
\eqref{eq:pi-pm} one sees that the spectrum of $\pi_\pm(b)$ is either
real positive or real negative, depending on the $\pm$~sign. Hence,
the algebras $\A(D^2_{q+})$ and $\A(D^2_{q-})$ describe the two
hemispheres of $\Sf^2_q$ and may be thought of as quantum disks, thus
justifying the notation $D_{q\pm}$.

There is a symbol map $\sg \: \A(D^2_{q\pm}) \to \A(\Sf^1)$ that maps
these ``noncommutative disks'' to their common boundary $\Sf^1$, which
is the equator of the equatorial Podle\'s sphere $\Sf^2_q$.
Explicitly, the symbol map is given as a $*$-homomorphism on the
generators of $\A(D_{q,\pm}^2)$ by
\begin{align}
\sg(r_\pm(a)) := u; \qquad \sg(r_\pm(b)) := 0,
\label{eq:symbol-map}
\end{align}
where $u$ is the unitary generator of $\A(\Sf^1)$.

Recall the algebra $\B$ defined around \eqref{eq:extended-pi} with
generators $\tilde a_\pm$, $\tilde b_\pm$ and
$P^\dn a_+ P^\up + P^\up a_- P^\dn$,
$P^\dn b_+ P^\up + P^\up b_- P^\dn$.
The following result emulates Proposition~4 of \cite{ConnesSUq} and establishes the correspondence \eqref{eq:corr-smoothing}. The results of~\cite{Naiad} on the approximate representation are crucial to its proof.

\begin{prop}
\label{pr:symbol-map}
There is a $*$-homomorphism
\begin{equation}
\rho: \B \to \A(D^2_{q+}) \ox \A(D^2_{q-}) \ox \A(\Sf^1)
\label{eq:symbol-map-bis}
\end{equation}
defined on generators by
\begin{align*}
\rho(\tilde a_+) &:= r_+(a) \ox r_-(a) \ox u,  &
\rho(\tilde a_-) &:= -q\,r_+(b) \ox r_-(b^*) \ox u^*,
\\
\rho(\tilde b_+) &:= - r_+(a) \ox r_-(b) \ox u, &
\rho(\tilde b_-) &:= - r_+(b) \ox r_-(a^*) \ox u^*.
\end{align*}
while the off-diagonal operators $P^\dn a_+ P^\up + P^\up a_- P^\dn$
and $P^\dn b_+ P^\up + P^\up b_- P^\dn$ are declared to lie in the
kernel of~$\rho$.
\end{prop}

\begin{proof}
First note that the $j$-dependence of the operators in $\B$ is
taken care of by the factor~$u$. Thus, it is enough to show that the
following prescription,
\begin{align*}
\rho_1(\tilde a_+) &:= \pi_+(a) \ox \pi_-(a), &
\rho_1(\tilde a_-) &:= -q\,\pi_+(b) \ox \pi_-(b^*),
\\
\rho_1(\tilde b_+) &:= - \pi_+(a) \ox \pi_-(b), &
\rho_1(\tilde b_-) &:= - \pi_+(b) \ox \pi_-(a^*),
\end{align*}
together with $\rho_1(P^\dn a_+ P^\up + P^\up a_- P^\dn) =
\rho_1(P^\dn b_+ P^\up + P^\up b_- P^\dn) := 0$, defines a
$*$-homo\-morph\-ism $\rho_1: \B \to \A(D^2_{q+}) \ox \A(D^2_{q-})$.
In the notation, we have replaced the representations $\pi_\pm$ of
$\A(SU_q(2))$ by corresponding faithful representations of
$\A(D^2_{q\pm})$ (omitting the maps $r_\pm$).

We define a map $\Pi: \H \to (\ell^2(\N) \otimes \ell^2(\N)) \otimes \C^2$, which simply forgets the $j$-index on the basis vectors $v_{xy}^j$: 
$$
\Pi: v_{xy}^j= \begin{pmatrix} v_{xy\up}^j \\ v_{xy\dn}^j \end{pmatrix} \mapsto \eps_{xy}:=\begin{pmatrix} \eps_{xy\up} \\ \eps_{xy\dn} \end{pmatrix},
$$ 
where $\eps_{xy\up}:=\eps_x \otimes \eps_y$ and $\eps_{xy\dn}:=\eps_x \otimes \eps_y$ in the two respective copies of $\ell^2(\N) \otimes \ell^2(\N)$ in its tensor product with  $\C^2$.

For any operator $T$ in $\B$, we define the map $\rho_1$ by
\begin{equation}
\rho_1(T) \eps_{xy} = \lim_{j \to \infty} \Pi(T v^j_{xy}).
\end{equation}
This map is well-defined, since $T$ is a polynomial in the generators of $\B$. Each such generator shifts the indices $x,y,j$ by $\pm \half$, with a coefficient matrix that can be bounded uniformly in $x,y$ and $j$ (cf. \cite{Naiad}) so that the limit $j \to \infty$ exists. 

First of all, it can be directly verified, using estimates
given in \cite[Sec.~7]{Naiad}, that the off-diagonal operators
$P^\dn a_+ P^\up + P^\up a_- P^\dn$ and $P^\dn b_+ P^\up + P^\up b_- P^
\dn$
are in the kernel of  $\rho_1$.
Next, the differences between the generators and the approximate 
generators $\ul a_\pm - \tilde a_\pm$
(and similarly $\tilde b_\pm - \ul b_\pm$)
lie in the kernel of $\rho_1$, as well.
Hence we can replace
$\tilde a_\pm$ and $\tilde b_\pm$ by $\ul a_\pm$
and $\ul b_\pm$, respectively.

Since the coefficients in the definition of $\ul a_\pm$ and $\ul b_\pm$ (equation \eqref{eq:approx-repn-bis}) are $j$-independent, we conclude that $\rho_1$ is of the desired form. For example, we compute:
\begin{align*}
\rho_1(\tilde a_+) \eps_{xy}=\rho_1(\ul a_+) \eps_{xy} &= \lim_{j \to \infty}  \sqrt{1-q^{2x+2}} \sqrt{1-q^{2y+2}} \Pi(v_{x+1, y+1}^{j^+})\\
&=  \sqrt{1-q^{2x+2}} \sqrt{1-q^{2y+2}} \eps_{x+1, y+1}= (\pi_+ (a) \otimes \pi_-(a) \otimes 1_2) \eps_{xy}.
\end{align*} 
Since a product of the operators $\ul a_\pm$ and $\ul b_\pm$ still does not contain $j$-dependent coefficients, $\rho_1$ respects the multiplication in $\B$. By linearity of the limit, $\rho_1$ is an algebra map. 
\end{proof}

\begin{defn}
The \textit{cosphere bundle on $SU_q(2)$} is defined as the range
of the map $\rho$ in $\A(D^2_{q+}) \ox A(D^2_{q-}) \ox \A(\Sf^1)$
and is denoted by $\A(\Sf_q^*)$.
\end{defn}

Note that $\Sf_q^*$ coincides with the cosphere bundle defined in
\cite{ConnesSUq,ConnesCIME}, where it is regarded as a noncommutative
space over which $D^2_{q+} \x D^2_{q-} \x \Sf^1$ is fibred.

The symbol map $\rho$ rectifies the correspondence
\eqref{eq:corr-smoothing}. Denote by $Q$ the orthogonal projector on
$\ell^2(\N)\ox \ell^2(\N) \ox \ell^2(\Z) \ox \C^2$ with range $\H'$,
which is the Hilbert subspace previously identified with~$\H$ just
before \eqref{eq:corr-smoothing}. Using \eqref{eq:corr-smoothing} in
combination with Proposition~\ref{pr:symbol-map}, we conclude that
\begin{equation}
T - Q (\rho(T) \ox 1_2) Q \in \OP^{-\infty}
\sepword{for all}  T \in \B.
\label{eq:symbol-smoothing}
\end{equation}
Here, the action of $\rho(T)$ on
$\ell^2(\N) \ox \ell^2(\N) \ox \ell^2(\Z)$ is determined by regarding
$\ell^2(\Z)$ as the Hilbert space of square-summable Fourier series
on~$\Sf^1$.

\section{The dimension spectrum}

We again follow \cite{ConnesSUq} for the computation of the dimension
spectrum. We define three linear functionals $\tau_0^\up$,
$\tau_0^\dn$ and $\tau_1$ on the algebras $\A(D_{q\pm}^2)$. Since
their definitions for both disks $D_{q+}^2$ and $D_{q-}^2$ are
identical, we shall omit the $\pm$ for notational convenience.

For $x \in \A(D^2_{q})$ we define,
\begin{align*}
\tau_1(x) &:= \frac{1}{2\pi} \int_{S^1} \sg(x),
\\
\tau_0^\up(x) &:= \lim_{N\to\infty} \Tr_N \pi(x) - (N+\sesq) \tau_1(x),
\\
\tau_0^\dn(x) &:= \lim_{N\to\infty} \Tr_N \pi(x) - (N+\half) \tau_1(x),
\end{align*}
where $\sg$ is the symbol map \eqref{eq:symbol-map}, and $\Tr_N$ is
the truncated trace
$$
\Tr_N(T) := \sum_{k=0}^N \braket{\eps_k}{T\eps_k}.
$$
The definition of the two different maps $\tau_0^\up$ and $\tau_0^\dn$
is suggested by the constants $\sesq$ and $\half$ appearing in our
choice of the Dirac operator; it will simplify some residue formulas
later on. We find that
\begin{align*}
\Tr_N(\pi(a))
&= (N + \sesq)\tau_1(a) + \tau_0^\up(a) + O(N^{-k})
\\
&= (N + \half)\tau_1(a) + \tau_0^\dn(a) + O(N^{-k})
\sepword{for all}  k > 0.
\end{align*}

Let us denote by $r$ the restriction homomorphism from
$\A(D^2_{q+}) \ox A(D^2_{q-}) \ox \A(\Sf^1)$ onto the first two
legs of the tensor product. In particular, we will use it as a map
$$
r: \A(\Sf_q^*) \to \A(D^2_{q+}) \ox A(D^2_{q-}).
$$
In the following, we adopt the notation~\cite{ConnesMIndex}:
$$
\ncint T := \Res_{z=0} \Tr T |D|^{-z}.
$$

\begin{thm}
The dimension spectrum of the spectral triple $(\A(SU_q(2)),\H,D)$ is
simple and given by $\{1,2,3\}$; the corresponding residues are
\begin{align*}
\ncint T |D|^{-3} &= 2(\tau_1 \ox \tau_1) \bigl(r\rho(T)^0\bigr),
\\
\ncint T |D|^{-2} &= \bigl(\tau_1 \ox (\tau_0^\up + \tau_0^\dn)
+ (\tau_0^\up + \tau_0^\dn) \ox \tau_1\bigr) \bigl(r\rho(T)^0\bigr),
\\
\ncint T |D|^{-1}
&= (\tau_0^\up \ox \tau_0^\dn + \tau_0^\dn \ox \tau_0^\up)
\bigl(r\rho(T)^0\bigr),
\end{align*}
with $T \in \Psi^0(\A)$.
\end{thm}

\begin{proof}
If we identify
$\H' \subset \ell^2(\N)\ox \ell^2(\N) \ox \ell^2(\Z) \ox \C^2$ with
$\H$ as above, the one-parameter group of automorphisms $\gamma(t)$
induces a $\Z$-grading on $\A(\Sf^*_q)$, in its representation on
$\H'$. We denote by $\rho(T)^0$ the degree-zero part of the diagonal
operator $\rho(T)$, for $T \in \B$. For the calculation of the
dimension spectrum we need to find the poles of the zeta function
$\zeta_T(z) := \Tr(T |D|^{-z})$ for all $T \in \Psi^0(\A)$. From our
discussion of the generators of $\Psi^0(\A)$, we see that we only need
to adjoin $P^\up \B$ to~$\B$.

In the zeta function $\zeta_T(z)$ for $T \in \B$, we can replace $T$
by $Q(\rho(T) \ox 1_2)Q$ since their difference is a smoothing
operator by \eqref{eq:symbol-smoothing}. The operator
$Q(\rho(T) \ox 1_2)Q$ commutes with the projector $P^\up$ so we can
first calculate
\begin{align}
\Tr(P^\up Q(\rho(T) \ox 1_2)Q\,|D|^{-z})
&= \sum_{2j=0}^\infty (2j + \sesq)^{-z}
(\Tr_{2j} \ox \Tr_{2j+1}) (r \rho(T)^0)
\nn \\
&= (\tau_1 \ox \tau_1) (r \rho(T)^0)\,\zeta(z - 2)
\nn \\
&\qquad + (\tau_1 \ox \tau_0^\dn + \tau_0^\up \ox \tau_1)
(r \rho(T)^0 )\,\zeta(z - 1)
\nn \\
&\qquad + (\tau_0^\up \ox \tau_0^\dn)
(r \rho(T)^0)\,\zeta(z) + f_\up(z),
\label{eq:res-compute-up}
\end{align}
where $f_\up(z)$ is holomorphic in $z \in \C$. Similarly,
\begin{align}
\Tr(P^\dn Q(\rho(T) \ox 1_2)Q\,|D|^{-z})
&= \sum_{2j=0}^\infty (2j + \sesq)^{-z}
(\Tr_{2j+1} \ox \Tr_{2j}) (r\rho(T)^0)
\nn \\
&= (\tau_1 \ox \tau_1) (r\rho(T)^0)\,\zeta(z - 2)
\nn \\
&\qquad + (\tau_1 \ox \tau_0^\up + \tau_0^\dn \ox \tau_1)
(r\rho(T)^0)\,\zeta(z - 1)
\nn \\
&\qquad + (\tau_0^\dn \ox \tau_0^\up)
(r\rho(T)^0)\,\zeta(z) + f_\dn(z),
\label{eq:res-compute-dn}
\end{align}
where $f_\dn(z)$ is holomorphic in~$z$. Since $\zeta(z)$ has a simple
pole at $z = 1$, we see that the zeta function $\zeta_T$ has simple
poles at $1$, $2$ and~$3$.
\end{proof}

{}From the above proof, we derive the following formulas which will be
used later on:
\begin{align}
\ncint P^\up T |D|^{-3} &= (\tau_1 \ox \tau_1) \bigl(r\rho(T)^0\bigr),
\nn \\
\ncint P^\up T |D|^{-2}
&= \bigl(\tau_1 \ox \tau_0^\dn + \tau_0^\up \ox \tau_1 \bigr)
\bigl(r\rho(T)^0\bigr),
\nn \\
\ncint P^\up T |D|^{-1}
&= (\tau_0^\up \ox \tau_0^\dn) \bigl(r\rho(T)^0\bigr),
\label{eq:expr-res-up}
\end{align}
with $T$ any element  in $\Psi^0(\A)$.

\section{Local index formula (\boldmath $d = 3$)}

We begin by discussing the local cyclic cocycles giving the local
index formula, in the general case when the spectral triple
$(\A,\H,D)$ has simple discrete dimension spectrum not containing $0$
and bounded above by~$3$.

Let us recall that with a general (odd) spectral triple $(\A,\H,D)$
there comes a Fredholm index of the operator $D$ as an additive map
$\varphi : K_1(\A) \to \Z$ defined as follows. If $F = \Sign D$ and
$P$ is the projector $P = \half(1 + F)$ then
\begin{equation}
\varphi([u])= \ind(PuP),
\label{eq:fr-ind}
\end{equation}
with $u \in \Mat_r(\A)$ a unitary representative of the $K_1$ class
(the operator $PuP$ is automatically Fredholm). The above map is
computed by pairing $K_1(\A)$ with ``nonlocal'' cyclic cocycles
$\chi_n$ given in terms of the operator $F$ and of the form
\begin{equation}
\chi_n(a_0, \dots, \a_n) = \lambda_n \Tr(a_0\,[F,a_1] \dots [F,a_n]),
\sepword{for all}  a_j \in \A,
\label{eq:nlcc}
\end{equation}
where $\lambda_n$ is a suitable normalization constant. The choice of
the integer $n$ is determined by the degree of summability of the
Fredholm module $(\H,F)$ over $\A$; any such module is declared to be
$p$-summable if the commutator $[F,a]$ is an element in the $p$-th
Schatten ideal $\L^p(\H)$, for any $a \in A$. The minimal $n$ in
\eqref{eq:nlcc} needs to be taken such that $n \geq p$.

On the other hand, the Connes--Moscovici local index theorem
\cite{ConnesMIndex} expresses the index map in terms of a local
cocycle $\phi_\odd$ in the $(b,B)$ bicomplex of $\A$ which is a local
representative of the cyclic cohomology class of $\chi_n$ (the cyclic
cohomology Chern character). The cocycle $\phi_\odd$ is given in terms
of the operator $D$ and is made of a finite number of terms
$\phi_\odd = (\phi_1, \phi_3, \dots )$; the pairing of the cyclic
cohomology class $[\phi_\odd] \in HC^\odd(\A)$ with $K_1(\A)$ gives
the Fredholm index \eqref{eq:fr-ind} of $D$ with coefficients in
$K_1(\A)$. The components of the cyclic cocycle $\phi_\odd$ are
explicitly given in \cite{ConnesMIndex}; we shall presently give them
for our case.

We know from Proposition~\ref{pr:regular} that our spectral triple
$(\A,\H,D)$ with $\A = \A(SU_q(2))$ has metric dimension equal to $3$.
As for the corresponding Fredholm module $(\H,F)$ over
$\A = \A(SU_q(2))$, it is $1$-summable since all commutators
$[F,\pi(x)]$, with $x\in\A$, are off-diagonal operators given by
sequences of rapid decay. Hence each $[F,\pi(x)]$ is trace-class and
we need only the first Chern character
$\chi_1(a_0,a_1) = \Tr(a_0\,[F,a_1])$, with $a_1,a_2 \in \A$ (we shall
omit discussing the normalization constant for the time being and come
back to it in the next section). An explicit expression for this
cyclic cocycle on the PBW-basis of $SU_q(2)$ was obtained
in~\cite{MasudaNW}.

The local cocycle has two components,
$\phi_\odd = (\phi_1, \phi_3)$, the cocycle condition
$(b + B)\phi_\odd = 0$ reading $B \phi_1=0,~ b \phi_1 + B \phi_3 =0,~b \phi_3= 0$ (see
Appendix~\ref{sec:cy-co}); it is explicitly given by
\begin{align*}
\phi_1(a_0,a_1) &:= \ncint a_0\, [D,a_1] \,|D|^{-1}
- \frac{1}{4} \ncint a_0\, \nabla([D,a_1]) \,|D|^{-3}
+ \frac{1}{8} \ncint a_0\, \nabla^2([D,a_1]) \,|D|^{-5},
\\
\phi_3(a_0,a_1,a_2,a_3)
&:= \frac{1}{12} \ncint a_0\,[D,a_1]\,[D,a_2]\,[D,a_3]\,|D|^{-3},
\end{align*}
where $\nabla(T) := [D^2,T]$ for any operator $T$ on~$\H$. Under the
assumption that $[F,a]$ is traceclass for each $a \in \A$, these
expressions can be rewritten as follows:
\begin{align}
\phi_1(a_0,a_1) &= \ncint a_0 \,\delta(a_1) F|D|^{-1}
- \frac{1}{2} \ncint a_0 \,\delta^2(a_1) F|D|^{-2}
+ \frac{1}{4} \ncint a_0 \,\delta^3(a_1) F|D|^{-3},
\nn \\
\phi_3(a_0,a_1,a_2,a_3) &= \frac{1}{12} \ncint a_0 \,\delta(a_1)
\,\delta(a_2) \,\delta(a_3) F|D|^{-3}.
\label{eq:odd-cycle}
\end{align}

We now quote Proposition~2 of \cite{ConnesSUq}, referring to that
paper for its proof.

\begin{prop}
\label{pr:LIF-one}
Let $(\A,\H,D)$ be a spectral triple with discrete simple dimension
spectrum not containing $0$ and bounded above by~$3$. If $[F,a]$ is
trace-class for all $a \in \A$, then the Chern character $\chi_1$ is
equal to $\phi_\odd - (b + B) \phi_\ev$ where the cochain
$\phi_\ev = (\phi_0, \phi_2)$ is given by
\begin{align*}
\phi_0(a) &:= \Tr(Fa\,|D|^{-z}) \bigr|_{z=0},
\\
\phi_2(a_0,a_1,a_2)
&:= \frac{1}{24} \ncint a_0 \,\delta(a_1) \,\delta^2(a_2) F|D|^{-3}.
\end{align*}
\end{prop}

The absence of $0$ in the dimension spectrum is needed for the
definition of $\phi_0$. The cochain
$\phi_\ev = (\phi_0, \phi_2)$ was named $\eta$-cochain in
\cite{ConnesSUq}. In components, the equivalence of the characters
means that
$$
\phi_1 = \chi_1 + b \phi_0 + B \phi_2 , \qquad
\phi_3 = b \phi_2 .
$$

The following general result, in combination with the above
proposition, shows that $\chi_1$ can be given (up to coboundaries) in
terms of one single $(b,B)$-cocycle~$\psi_1$.

\begin{prop}
\label{pr:LIF-two}
Let $(\A,\H,D)$ be a spectral triple with discrete simple dimension
spectrum not containing $0$ and bounded above by~$3$. Assume that
$[F,a]$ is trace class for all $a \in \A$, and set
$P := \half(1 + F)$. Then, the local Chern character
$\phi_\odd$ is equal to $\psi_1 - (b + B)\phi'_\ev$, where
$$
\psi_1(a_0,a_1) := 2 \ncint a_0 \,\delta(a_1) P|D|^{-1}
- \ncint a_0 \,\delta^2(a_1) P|D|^{-2}
+ \frac{2}{3} \ncint a_0 \,\delta^3(a_1) P|D|^{-3},
$$
and $\phi'_\ev = (\phi'_0, \phi'_2)$ is given by
\begin{align*}
\phi'_0(a) &:= \Tr(a\,|D|^{-z}) \bigr|_{z=0},
\\
\phi'_2(a_0, a_1, a_2)
&:= - \frac{1}{24} \ncint a_0 \,\delta(a_1) \,\delta^2(a_2) F|D|^{-3}.
\end{align*}
\end{prop}

\begin{proof}
One needs to verify the following equalities between cochains in the
$(b,B)$ bicomplex:
\begin{align*}
\phi_1 + b \phi'_0 + B \phi'_2 &= \psi_1,
\\
\phi_3 + b \phi'_2 &= 0.
\end{align*}
The second equality follows from a direct computation of $b\phi'_2$
and comparing with equation \eqref{eq:odd-cycle}. Note that this
identity proves that $\psi_1$ is indeed a cyclic cocycle. One also
shows that
$$
B \phi'_2(a_0,a_1) = \frac{1}{12} \ncint a_0\,\delta^3(a_1) F|D|^{-3}.
$$
Then, using the asymptotic expansion \cite{ConnesMIndex}:
$$
|D|^{-z} a \sim \sum_{k\geq 0} \binom{-z}{k} \delta^k(a)\,|D|^{-z-k}
$$
modulo very low powers of $|D|$, one computes
$$
b \phi'_0(a_0,a_1) = \ncint a_0\,\delta(a_1) |D|^{-1}
- \frac{1}{2} \ncint a_0 \,\delta^2(a_1) |D|^{-2}
+ \frac{1}{3} \ncint a_0 \,\delta^3(a_1) |D|^{-3},
$$
and it is now immediate that $\phi_1 + b \phi'_0 + B \phi'_2$ gives
the cyclic cocycle $\psi_1$.
\end{proof}

\begin{rem}
The term involving $P |D|^{-3}$ would vanish if the latter were
traceclass, which is the case in \cite{ConnesSUq} (this is the
statement that the metric dimension of the projector $P$ is~$2$).
\end{rem}

Combining these two propositions, it follows that
the cyclic $1$-cocycles $\chi_1$ and $\psi_1$ are related as:
\begin{equation}
\chi_1 = \psi_1 -  b \beta,
\label{eq:tau-psi}
\end{equation}
where $\beta (a) = 2\Tr(Pa\,|D|^{-z})\bigr|_{z=0}$.

\section{The pairing between \boldmath $HC^1$ and $K_1$}

In this section, we shall calculate the value of the index map
\eqref{eq:fr-ind} when $U$ is the unitary operator representing the
generator of $K_1(\A(SU_q(2)))$,
$$
\varphi([U]) = \ind(PUP) := \dim\ker P UP - \dim\ker P U^* P,
$$
with
\begin{equation}
U = \begin{pmatrix} a & b \\ -q b^* & a^* \end{pmatrix},
\label{eq:uni-def}
\end{equation}
acting on the doubled Hilbert space $\H \ox \C^2$ via the
representation $\pi \ox 1_2$. The projector $P$ was denoted $P^\up$ in
Section~\ref{sec:iso-ge}. One expects this index to be nonzero, since the $K$-homology class of $(\A,\H,D)$ is non-trivial. This has been remarked also in \cite{ChakrabortyPDirac}, where our spectral triple is decomposed in terms of the spectral triple constructed in \cite{ChakrabortyPEqvt}. 

We first compute the above index directly, which is possible due to the simple nature of this particular example.
A short computation shows that the kernel of the operator $P U^* P$ is
trivial, whereas the kernel of $P U P$ contains only elements
proportional to the vector
$$
\begin{pmatrix} \quad \ket{0,0,-\half, \up} \\
 -q^{-1} \ket{0,0,\half,\up}
\end{pmatrix},
$$
leading to $\varphi([U]) = \ind(PUP) = 1$.

\vspace{6pt}

Recall that for $\A = \A(SU_q(2))$, our Fredholm module $(\H,F)$ over
$\A(SU_q(2))$ is $1$-summable. From the previous section we know that
$\ind(PUP)$ can be computed using the local cyclic cocycle
$\psi_1$, see eqn.~\eqref{eq:tau-psi}. To prepare for this index
computation via $\psi_1$, we recall the following lemma
\cite[IV.1.$\ga$]{Book}, which fixes the normalization constant in
front of $\chi_1$. For completeness we recall the proof.

\begin{lem}
Let $(\H,F)$ be a $1$-summable Fredholm module over $\A$ with
$P = \half(1 + F)$; let $u  \in \Mat_r(\A)$ be unitary with a suitable
$r$. Then $P u P$ is a Fredholm operator on $P \H$ and
$$
\ind(P u P) = -\half \Tr(u^* [F,u]) = -\half \chi_1(u^*,u).
$$
\end{lem}

\begin{proof}
We claim that $P u^* P$ is a parametrix for $P u P$, that is, an
inverse modulo compact operators on $P\H$. Indeed, since
$P - u^* P u = -\half u^*\,[F,u]$ is traceclass by 1-summability, by
composing it from both sides with $P$ it follows that $P - P u^* P u P$
is traceclass. Therefore,
\begin{equation}
\ind(P u P) = \Tr(P - P u^* P u P) - \Tr(P - P u P u^* P),
\label{eq:index-formula}
\end{equation}
and the identities $P - P u^* P u P = -\half P u^*\,[F,u] P$ and
$[F,u]\,u^* + u\,[F,u^*] = 0$, together with $[F,[F,u]]_+=0$, imply the statement.
\end{proof}

Thus, the index of $P U P$, for the $U$ of~\eqref{eq:uni-def} is
given, up to an overall $-\half$ factor, by
$$
\psi_1(U^{-1}, U)
= 2 \ncint U_{kl}^* \,\delta(U_{lk}) P |D|^{-1}
- \ncint U_{kl}^* \,\delta^2(U_{lk}) P |D|^{-2}
+ \frac{2}{3} \ncint U_{kl}^* \,\delta^3(U_{lk}) P |D|^{-3},
$$
with summation over $k,l = 0,1$ understood. We compute this expression
using equation \eqref{eq:expr-res-up}. First note that since the
entries of $U$ are generators of $\A(SU_q(2))$, we see from
\eqref{eq:delta-apm} and~\eqref{eq:delta-bpm} that
$\rho(\delta^2(U_{kl})) = \rho(U_{kl})$, a relation that simplifies
the above formula. We compute the degree~$0$ part of
$\rho(U_{kl}^*\,\delta(U_{lk}))$ with respect to the grading coming
from $\ga(t)$ --the only part that contributes to the trace-- using
the algebra relations of $\A(D_{q\pm}^2)$,
$$
\rho(U_{kl}^*\,\delta(U_{lk}))^0 = 2(1 - q^2)\,1 \ox r_-(b)^2.
$$
Using the basic equalities
$$
\tau_1(1) = 1,  \quad \tau_1(r_\pm(b)^n) = 0, \quad
\tau_0^\up(1) = -\tau_0^\dn(1) = -\half,  \quad
\tau_0^\up(r_\pm(b)^n) = \tau_0^\dn(r_\pm(b)^n)
= \frac{(\pm 1)^n}{1 - q^n},
$$
we find that
$$
\psi_1(U^{-1},U)
= 2(1 - q^2) (2\tau_0^\up \ox \tau_0^\dn + \frac{2}{3}\tau_1 \ox \tau_1)
\bigl(1 \ox r_-(b)^2 \bigr)
- (\tau_1 \ox \tau_0^\dn + \tau_0^\up \ox \tau_1)
\bigl(1 \ox 1 \bigr) = -2.
$$
Taking the proper coefficients, we finally obtain
$$
\ind(PUP) = -\half \psi_1(U^{-1}, U) = 1.
$$

\appendix
\section{Pseudodifferential calculus and cyclic cohomology}
\label{sec:pdc}
\label{sec:cy-co}

Recall \cite{CareyPRS,ConnesMIndex,Polaris} that a spectral triple
$(\A,\H,D)$ is \textit{regular} (or \textit{smooth}, or $QC^\infty$)
if the algebra generated by $\A$ and $[D,\A]$ lies within the smooth
domain $\bigcap_{n=0}^\infty \Dom \delta^n$ of the operator derivation
$\delta(T) := |D|T - T|D|$. This condition permits to introduce the
analogue of Sobolev spaces $\H^s := \Dom(1 + D^2)^{s/2}$ for
$s \in \R$. Let $\H^\infty := \bigcap_{s\geq 0} H^s$, which is a core
for~$|D|$. Then $T\: \H^\infty \to \H^\infty$ has \textit{analytic
order $\leq k$} if $T$ extends to a bounded operator from $\H^{k+s}$ to
$\H^s$ for all $s \geq 0$. It turns out that
$\A(\H^\infty) \subset \H^\infty$.

Assume that $|D|$ is invertible --which is a generic case of the $D$
used in this paper (for a careful treatment of the noninvertible case,
see \cite{CareyPRS}). The space $\OP^\a$ of \textit{operators of
order}~$\leq\a$ consists of those $T \: \H^\infty \to \H^\infty$ such
that
$$
|D|^{-\a} T \in \bigcap_{n=1}^\infty \Dom \delta^n.
$$
(Operators of order $\a$ have analytic order $\a$). In particular,
$\OP^0 = \bigcap_{n=1}^\infty \Dom \delta^n$, the algebra of operators
of order~$\leq 0$ includes $A \cup [D,\A]$ and their iterated
commutators with~$|D|$. Moreover, $[D^2, \OP^\a] \subset \OP^{\a+1}$
and $\OP^{-\infty} := \bigcap_{\a\leq 0} \OP^\a$ is a two-sided ideal
in~$\OP^0$.

The algebra structure can be read off in terms of an
\textit{asymptotic expansion}: $T \sim \sum_{j=0}^\infty T_j$
whenever $T$ and each $T_j$ are operators from $\H^\infty$ to
$\H^\infty$; and for each $m \in \Z$, there exists $N$ such that for
all $M > N$, the operator $T - \sum_{j=1}^M T_j$ has analytic order
$\leq m$. For instance, for complex powers of $|D|$ (defined by the
Cauchy formula) there is a binomial expansion:
$$
[|D|^z, T] \sim \sum_{k=1}^\infty \binom{z}{k} \delta^k(T)\,|D|^{z-k} .
$$

Thus far, we have employed finitely generated algebras $\A(X)$, where
$X = SU_q(2)$, $D^2_{q\pm}$, $\Sf^1$ or~$\Sf_q^2$. In each case, we
can enlarge them to algebras $C^\infty(X)$ by replacing polynomials
in the generators (given in a prescribed order) by series with
coefficients of rapid decay: this is clear when $X = \Sf^1$, where
smooth functions have rapidly decaying Fourier series. Using the
symbol maps \eqref{eq:q-disks}, \eqref{eq:symbol-map}
and \eqref{eq:symbol-map-bis} together with Lemma~2 of
\cite{ConnesCIME}, we can check that each such $C^\infty(X)$ is closed
under holomorphic functional calculus. The foregoing results apply,
mutatis mutandis, to the regular spectral triple
$(C^\infty(SU_q(2)),\H,D)$.

\vspace{8pt}

For convenience, we also summarize here the cyclic cohomology of the algebra $\A(SU_q(2))$. A cyclic $n$-cochain on an algebra $\A$ is an element
$\varphi \in \chala{n}$, the collection of $(n+1)$-linear
functionals on $\A$ which in addition are cyclic,
$\lambda \varphi = \varphi$, with 
$$
\lambda\varphi(a_0, a_1,\dots, a_n)
= (-1)^n \varphi(a_n, a_0,\dots, a_{n-1}).
$$
There is a cochain complex
$(C_\la^\8(\A) = \bigoplus_n \chala{n},\,b)$ with (Hochschild)
coboundary operator $b \: \cha{n} \to \cha{n+1}$ defined by
$$
b\varphi(a_0, a_1,\dots, a_{n+1})
:= \sum_{j=0}^n (-1)^j \varphi(a_0,\dots, a_j a_{j+1},\dots, a_{n+1})
+ (-1)^{n+1} \varphi(a_{n+1} a_0, a_1,\dots, a_n) .
$$
The cyclic cohomology $HC^\8(\A)$ of the algebra $\A$ is the
cohomology of this complex,
$$
HC^n(\A) :=  H^n(C_\la^\8(\A),b).
$$
Equivalently, $HC^\8(\A)$ can be described \cite{Book,Polaris} by
using the second filtration of a $(b,B)$ bicomplex of arbitrary (i.e.,
noncyclic) cochains on $\A$. Here the operator $B$ decreases the
degree $B\: \cha{n} \to \cha{n-1}$, and is defined as $B = N B_0$,
with
\begin{align*}
& (B_0 \varphi)(a_0,\dots, a_{n-1})
:= \varphi(1,a_0,\dots, a_{n-1}) - (-1)^n \varphi(a_0,\dots,a_{n-1},1)
\\
 & (N\psi)(a_0,\dots, a_{n-1}) := \sum_{j=0}^{n-1}
    (-1)^{(n-1)j} \psi(a_j,\dots,a_{n-1}, a_0,\dots,a_{j-1}).
\end{align*}
It is straightforward to check that $B^2 = 0$ and that $b B + B b = 0$;
thus $(b + B)^2 = 0$. By putting together these two operators, one
gets a bicomplex $(C^\8(\A), b, B)$ with $\cha{p-q}$ in bidegree
$(p,q)$. To a cyclic $n$-cocycle one associates the $(b,B)$ cocycle
$\varphi$, $(b + B)\varphi = 0$, having only one nonvanishing
component $\varphi_{n,0}$ given by
$\varphi_{n,0} := (-1)^{\piso{n/2}} \psi$.

\vspace{6pt}

The cyclic cohomology of the algebra $\A(SU_q(2))$ was computed in
\cite{MasudaNW}. The even components vanish while the odd ones were
found to be one-dimensional and generated by the cyclic $1$-cocycle
$\tau_\odd \in HC^1(\A(SU_q(2)))$ which was obtained as a character of
a $1$-summable Fredholm module,
$$
\tau_\odd(a^l b^m (b^*)^n, a^{l'} b^{m'} (b^*)^{n'}) = (n-m) \dfrac{q^{l(m'+n')} \prod_{i=1}^l (1 - q^{2i})}
        {\prod_{i=0}^l (1- q^{2i+2n+2n'})} \,\delta_{n+n', m+m'} \delta_{l,-l'}
$$
where we use the notation $a^{-l} = (a^*)^l$ for $l > 0$. Since $HC^1(\A(SU_q(2)))$ is one-dimensional, the characters of the 1-summable Fredholm modules found in \cite{ConnesSUq} and in this paper, are all cohomologous to this cyclic cocycle.

\vspace{6pt}

\subsection*{Acknowledgments}

We thank Alain Connes for suggesting that we address this problem. GL
and LD acknowledge partial support by the Italian Co-Fin Project
SINTESI. Support from the Vicerrector\'{\i}a de Investigaci\'on of the
Universidad de Costa Rica is acknowledged by JCV.


\begin{thebibliography}{11}

\bibitem{CareyPRS}
A. L. Carey, J. Phillips, A. Rennie and F. A. Sukochev,
``The Hochschild class of the Chern character for semifinite spectral
triples'',
J. Funct. Anal. {\bf 213} (2004), 111--153.

\bibitem{ChakrabortyPEqvt}
P. S. Chakraborty, A. Pal,
``Equivariant spectral triples on the quantum $SU(2)$ group'',
K-Theory {\bf 28} (2003), 107--126.

\bibitem{ChakrabortyPDirac}
P. S. Chakraborty, A. Pal,
``On equivariant Dirac operators for $SU_q(2)$'',
math.qa\slash 0501019, New Delhi, 2005.

\bibitem{Book}
A. Connes,
\textit{Noncommutative Geometry},
Academic Press, London and San Diego, 1994.

\bibitem{ConnesSUq}
A. Connes,
``Cyclic cohomology, quantum group symmetries and the local index
formula for $SU_q(2)$'',
J. Inst. Math. Jussieu {\bf 3} (2004), 17--68.

\bibitem{ConnesCIME}
A. Connes,
``Cyclic cohomology, noncommutative geometry and quantum group
symmetries'',
in \textit{Noncommutative Geometry},
S. Doplicher and R. Longo, eds.,
Lecture Notes in Mathematics {\bf 1831},
Springer, Berlin, 2004.

\bibitem{ConnesMIndex}
A. Connes and H. Moscovici,
``The local index formula in noncommutative geometry'',
Geom. Funct. Anal. {\bf 5} (1995), 174--243.

\bibitem{Naiad}
L.~D\c{a}browski, G.~Landi, A.~Sitarz, W.~van Suijlekom, and J.~C. V\'arilly. 
``The Dirac operator on $SU_q(2)$'',
Commun. Math. Phys. {\bf 259} (2005), 729-759.

\bibitem{Polaris}
J. M. Gracia-Bond\'{\i}a, J. C. V\'arilly and H. Figueroa,
\textit{Elements of Noncommutative Geometry},
Birkh\"auser, Boston, 2001.

\bibitem{MasudaNW}
T. Masuda, Y. Nakagami and J. Watanabe,
``Noncommutative differential geometry on the quantum $SU(2)$.
I: An algebraic viewpoint'',
K-Theory {\bf 4} (1990), 157--180.

\bibitem{Podles}
P. Podle\'s,
``Quantum spheres'',
Lett. Math. Phys. {\bf 14} (1987), 521--531.

\end{thebibliography}
\end{document}